\documentclass[12pt]{article}
\usepackage{amssymb, amsmath, amsfonts}
\usepackage{bm}
\usepackage[all,2cell]{xy}\UseAllTwocells
\long\def\nodo#1{}
\def\MM{{\mathcal M}}
\def\NN{{\mathcal N}}
\def\PP{{\mathcal P}}
\def\RR{{\mathcal R}}
\def\CC{{\mathcal C}}
\def\DD{{\mathcal D}}
\def\tto{{\Rightarrow}}
\def\Id{\mathrm{Id}}
\def\id{\mathrm{id}}
\def\Ob{\operatorname{Ob}}
\def\actu{\mathrm{\bf -act}^c_1}
\def\act{\mathrm{\bf -act}^c}
\def\actp{\mathrm{\bf -act}^p}
\def\cat{\mathrm{\bf Cat}}
\def\End{\operatorname{End}}
\def\uEnd{\operatorname{\bf End}}

\def\distr{\mathfrak{distr}}
\def\nxpoint{\vskip .1in\refstepcounter{section}%
{\bf\thesection. }}
\def\refpoint#1{{\rm\textbf{\ref{#1}}}}

\begin{document}
\begin{center}
{\it \large Equivariant monads and equivariant lifts versus
a 2-category of distributive laws}
\\
{\sc Zoran \v{S}koda} {(preliminary notes)}

\end{center}

\footnotesize{\bf
Fix a monoidal category $\CC$. The 2-category of monads in
the 2-category of $\CC$-actegories, colax $\CC$-equivarant functors,
and $\CC$-equivariant natural transformations of colax functors,
may be recast in terms of pairs consisting of
a usual monad and a distributive law between the monad and
the action of $\CC$, morphisms of monads respecting the distributive
law, and transformations of monads satisfying some compatibility with
the actions and distributive laws involved. The monads in this picture
may be generalized to actions of monoidal categories, and actions
of PRO-s in particular. If $\CC$ is a PRO as well, then in special
cases one gets various distributive laws of  a given classical type,
for example between a comonad and an endofunctor or 
between a monad and a comonad. The usual pentagons are in general
replaced by multigons, and there are also ``mixed'' multigons involving
two distinct distributive laws.
Beck's bijection between the distributive laws
and lifts of one monad to the Eilenberg-Moore category of another
monad is here extended to an isomorphism of 2-categories. The lifts
of maps of above mentioned pairs are colax $\CC$-equivariant.
We finish with a short treatment of relative distributive laws 
between two pseudoalgebra structures which are relative with 
respect to the distributivity of two pseudomonads involved, what gives
a hint toward the generalizations. 
}
\vskip .07in

\nxpoint Throughout the paper, $\CC$ will be a fixed monoidal category
with a monoidal product $\otimes$, a unit object $\bm 1$, the associativity
coherence isomorphisms $a_{X,Y,Z} : X\otimes (Y\otimes Z)
\to (X\otimes Y)\otimes Z$,
natural in $X,Y,Z\in\mathrm{Ob}\,\CC$, the left unit coherence
$r : {\rm Id}_\CC \tto {\rm Id}_\CC\otimes {\bf 1}$
and the right unit coherence
$l : {\rm Id}_\CC \tto {\bf 1} \otimes {\rm Id}_\CC$
satisfying for all $A,B,C,D\in\Ob\CC$ the MacLane pentagon 
$a_{A,B,C\otimes D}\circ a_{A\otimes B,C,D}\circ (a_{A,B,C}\otimes D) 
= (A\otimes a_{B,C,D})\circ a_{A,B\otimes C,D}$
and unit triangle coherence relations
$a_{\bm 1, A,B}\circ l_{A\otimes B} = l_A\otimes B$
and $r_{A\otimes B} = a_{A,B,\bm 1}\circ (A\otimes r_B)$.
A left coherent action of $\CC$ on a category $\NN$
is a coherent monoidal functor
${\mathcal L} : \CC\to \End\NN$ where $\End\CC$ is strict 
monoidal with respect to the composition of endofunctors. Equivalently,
a $\CC$-action will be given by a bifunctor $\lozenge : \CC\times\NN\to\NN$,
natural isomorphisms
$\Psi : (\_\otimes\_)\lozenge\_\tto\_\lozenge(\_\lozenge\_)$
and $u : \Id_\NN\tto \bm 1\lozenge \Id_\NN$ satisfying for all 
$Q,Q',Q''\in\Ob\CC$ and $N\in\Ob\NN$ the
action pentagon coherence $\Psi_{Q,Q',Q''\lozenge N}\circ\Psi_{Q\otimes
Q', Q'',N}\circ (a_{Q,Q',Q''}\lozenge N) = (Q\lozenge \Psi_{Q',Q'',N})
\circ\Psi_{Q,Q'\otimes Q'',N}$ and unit action coherences
$u_{Q\lozenge N}\circ\Psi_{\bm 1,Q,N} = l_Q\lozenge N$
and $(Q\lozenge u_N)\circ\Psi_{Q,\bm 1,N}=r_Q\lozenge\bm N$.
A $\CC$-{\bf actegory} is a category $\NN$ equipped 
with a coherent action $\lozenge,\Psi, u$ of $\CC$.

\nxpoint A {\bf colax $\CC$-equivariant functor of $\CC$-actegories}
$(F,\zeta) : (\MM,\lozenge^\MM,\Psi^\MM,u^\MM) 
\to (\NN,\lozenge^\NN,\Psi^\NN,u^\NN)$  
is a usual functor $F :\MM\to\NN$ 
with a binatural transformation of bifunctors
$\zeta : F(\_\lozenge^\MM\_)\tto \_\tilde\lozenge^\NN F(\_)
:\CC\times\MM\tto\NN$, so that 
\begin{equation}\label{eq:colaxu}
\xymatrix{
\bm F(\bm 1\lozenge^\MM M) \ar[rr]^{\zeta_{\bm 1,M}} 
\ar[rd]_{F(u^\MM_M)}
&& \bm 1\lozenge^\NN F(M)\ar[ld]^{u^\NN_{F(M)}}\\
& F(M) &
}\end{equation}
\begin{equation}\label{eq:colaxPsi}\xymatrix{
F((A\otimes B)\lozenge^\MM M)\ar[rr]^{\zeta_{(A\otimes B)\lozenge M}}
\ar[d]_{F(\Psi^\MM_{A,B,M})}
&& (A\otimes B)\lozenge^\NN F(M)\ar[d]^{\Psi^\NN_{A,B,F(M)}}
\\
F(A\lozenge^\MM (B\lozenge^\NN M))\ar[r]^{\zeta_{A,B\lozenge M}} &
A\lozenge^\NN F(B\lozenge^\MM M) \ar[r]^{A\lozenge \zeta_{B,N}} &
A\lozenge^\NN (B\lozenge^\NN F(M))
} 
\end{equation}
$\CC$-actegories and colax $\CC$-equivariant functors make a category
$\CC\actu$:  
given $(F,\zeta^F) : \NN\to\PP$ and $(G,\zeta^G):\MM\to \NN$ their
composition is
$(F,\zeta^F)\circ(G,\zeta^G) := (F\circ G, \zeta^{F\circ G}):\MM\to\PP$ 
where $\zeta^{F\circ G}_{C,M} := \zeta^F_{C,GM}\circ\, F(\zeta_{C,M}) 
: F(G(C\lozenge^\MM M))\tto C\lozenge^\PP F(G(M))$.

\nxpoint 
A {\bf $\CC$-equivariant natural transformation of colax $\CC$-equivariant
functors}
$\alpha : (F,\zeta^F)\tto (H,\zeta^H) : \MM\to\NN$
is a natural transformation of underlying ordinary functors 
$\alpha : F\tto G$ such that for all
$C\in \CC, M\in M$ the following square commutes:
\begin{equation}\label{eq:equivCCnat}\xymatrix{
F(C\lozenge^\MM M) \ar[r]^{\zeta^F_{C,M}} \ar[d]_{\alpha_{C\lozenge M}}
& C\lozenge^\NN FM \ar[d]^{C\lozenge\alpha_M}\\
G(C\lozenge^\MM M) \ar[r]^{\zeta^G_{C,M}} 
& C\lozenge^\NN GM
}\end{equation}
The usual transformation of usual functors obtained as 
a vertical or a horizontal composition of 
$\CC$-equivariant natural transformations 
of colax $\CC$-functors is $\CC$-equivariant. 
Thus we obtain a strict 2-category $\CC\act$ which has
all cartesian products, namely the usual products in $\mathrm{\bf Cat}$
equipped with the diagonal $\CC$-action, e.g. for binary products
of $\CC$-actegories $C\lozenge (M,N) = (C\lozenge M,C\lozenge N)$,
and for $\CC$-functors 
$(F,\zeta^F)\times(G,\zeta^G)=(F\times G,\zeta^F\times\zeta^G)$.
\nxpoint Let $G$ be an endofunctor on a category $\MM$. For a given monad
$\bf T = (T,\mu,\eta)$ with a multiplication $\mu : TT\tto T$
and unit $\nu : \Id\tto T$ a {\bf distributive law} between $G$ and $T$
is a natural transformation $l : GT\tto TG$ such that
\begin{equation}\tag{D1}
\xymatrix{
GTT \ar[d]^{G\mu}\ar[r]^{lT} & TGT \ar[r]^{Tl} & TTG \ar[d]^{\mu G}\\
GT\ar[rr]^{l} && TG
}\end{equation}
commutes and $l \circ G\eta = \eta G : G \tto TG$.
A {\bf lift} of an endofunctor (resp. (co)monad) $G$
to a category $\mathrm C$ equipped with a functor $U$ to
$\MM$ is an endofunctor (resp. (co)monad) $\tilde{G}$
such that $U \tilde{G} = GU$ (and obvious additonal conditions
for the (co)monad case).
The basic motivating fact for this definition states
that the distributive laws between $G$ and $T$
are in a canonical bijection with the lifts of
endofunctor $G$ to the Eilenberg-Moore category
$\MM^{\bf T}$ of modules $(M,\nu)$
with respect to the forgetful functor $U : (M,\nu)\mapsto M$
(as usual, $M \in\Ob\MM$ and $\nu : TM\to M$).
Often $G$ is also a (co)monad.
Then, two additional axioms are required for $l$
which ensure that $\tilde{G}$ is also a (co)monad.
Modulo quoting this very fact, no proof in this paper
needs repair when replacing distributive laws and lifts where
$G$ is endofunctor, with the version where $G$ is a (co)monad.
\nxpoint In every strict 2-category, endo-1-cells of a fixed object
and their natural transformations form a strict monoidal category,
with the horizontal composition as the tensor product.
In particular, $\uEnd_\CC(\MM):=\CC\act(\MM,\MM)$ is
a strict monoidal category. If $\tilde{T} = (T, \zeta)$ is an object 
in $\uEnd(\MM)$, that is a colax $\CC$-equivariant endofunctor, then
its tensor square is 
$\tilde{T}\tilde{T} := (T\circ T, \zeta_T \circ T(\zeta))$.
Here $(\zeta_T\circ T(\zeta))_{C,M} := \zeta_{C,TM}\circ T(\zeta_{C,M})
: TT(C\lozenge M)\tto C\lozenge TTM$. 
Let now $\tilde{\bf T} = (\tilde{T},\mu,\eta)$
be a monad in $\CC\act$. Our next aim is to decipher these data in terms
of data in $\mathrm{\bf Cat}$. 
$\tilde{T} = (T,\zeta)$ is a colax $\CC$-equivariant endofunctor hence
the two diagrams~(\ref{eq:colaxu},\ref{eq:colaxPsi}) commute
with $T$ in place of $F$. The multiplication 
$\mu : \tilde{T}\tilde{T}\tto\tilde{T}$ is a natural transformation
$\mu :TT\tto T$, whose $\CC$-equivariance says 
that~(\ref{eq:equivCCnat}) commutes for
$\alpha = \mu$, $F = TT$, $G = T$, $\zeta^F = \zeta_T \circ T(\zeta)$
and $\zeta^G = \zeta$. From this we obtain the following pentagon
\begin{equation}\label{eq:TTCdistr}\xymatrix{
TT(C\lozenge M) \ar[r]^{T(\zeta_{C,M})}\ar[d]_{\mu_{C\lozenge M}} & 
T(C\lozenge TM)\ar[r]^{\zeta_{C,TM}} &
C\lozenge TTM\ar[d]^{C\lozenge \mu_M}
\\ T(C\lozenge M) \ar[rr]^{\zeta_{C,M}}&& C\lozenge TM
}\end{equation}
The unit $\eta :(\Id_\MM,\Id_\Id)\tto(T,\zeta)$ is a natural tranformation
$\eta :\Id\to T$ and its $\CC$-equivariance means that~(\ref{eq:equivCCnat})
commutes for $F=\Id_\MM$, $\zeta^F=\Id_\Id$, $G=T$ and $\zeta^G=\zeta$
what reduces to the triangle
\begin{equation}\label{eq:etadistr}\xymatrix{
& C\lozenge M\ar[ld]_{\eta_{C\lozenge M}}\ar[rd]^{T(\eta_M)}&\\
T(C\lozenge M)\ar[rr]^{\zeta_{C,M}}&& C\lozenge TM
}\end{equation}
The identities for $\mu$ and $\eta$ 
(monad associativity  $\mu \circ T(\mu) =\mu \circ \mu_T$and unit axioms)
simply say that the underlying endofunctor has a structure of a monad.

{\bf Proposition.} {\it A monad $\tilde{\bf T} = (\tilde{T},\mu,\eta)$
in $\CC\act$ is the same as a usual monad ${\bf T} = (T,\mu,\eta)$
together with a binatural transformation $\zeta : T(\_\lozenge\_)
\tto \_\lozenge T(\_)$, satisfying~(\ref{eq:colaxu}),(\ref{eq:colaxPsi}) with
$T=F$ and~(\ref{eq:TTCdistr}),(\ref{eq:etadistr}), i.e. the {\bf distributive
law} between $\CC$-action and ${\bf T}$. 
}

\nxpoint More generally, we may be given two actions of monoidal
categories $\CC$ and $\DD$ on the same category $\MM$. The distributive
law between these two actions will be a binatural transformation
of bifunctors $\DD\lozenge(\CC\lozenge\MM)\tto\CC\lozenge(\DD\lozenge\MM)$
satisfying again some coherences; this general case will be studied 
elsewhere and, in the case of one left and one right action 
also in~\cite{skoda:biact}. Let us now recall the classical
case. 

Let $G$ be an endofunctor on a category $\MM$. For a given monad
$\bf T = (T,\mu,\eta)$ with a multiplication $\mu : TT\tto T$
and unit $\nu : \Id\tto T$ a {\bf distributive law} between $G$ and $T$
is a natural transformation $l : GT\tto TG$ such that
\begin{equation}\tag{D1}
\xymatrix{
GTT \ar[d]^{G\mu}\ar[r]^{lT} & TGT \ar[r]^{Tl} & TTG \ar[d]^{\mu G}\\
GT\ar[rr]^{l} && TG
}\end{equation}
commutes and $l \circ G\eta = \eta G : G \tto TG$.
A {\bf lift} of an endofunctor (resp. (co)monad) $G$
to a category $\mathrm C$ equipped with a functor $U$ to
$\MM$ is an endofunctor (resp. (co)monad) $\tilde{G}$
such that $U \tilde{G} = GU$ (and obvious additonal conditions
for the (co)monad case).
The basic motivating fact for this definition states
that the distributive laws between $G$ and $T$
are in a canonical bijection with the lifts of
endofunctor $G$ to the Eilenberg-Moore category
$\MM^{\bf T}$ of modules $(M,\nu)$
with respect to the forgetful functor $U : (M,\nu)\mapsto M$
(as usual, $M \in\Ob\MM$ and $\nu : TM\to M$).
Often $G$ is also a (co)monad.
Then, two additional axioms are required for $l$
which ensure that $\tilde{G}$ is also a (co)monad.
The generalizations of these additional axioms for the 
case of PRO are also studied below. We start with the easier
partial case of monads.

\nxpoint \label{sec:alphaHcorr}
 A {\bf map of monads} in a fixed category $\MM$
is a natural transformation
$\alpha : T \Rightarrow T'$ for which $\alpha\circ \mu_T = \mu_{T'}$
and $\mu \circ T\eta = \mu\circ \eta T = \id : T \Rightarrow T$.
Every map of monads $\alpha$
induces a functor of Eilenberg-Moore categories
$H^\alpha : {\mathcal M}^{\bf T'}\to{\mathcal M}^{\bf T}$
by the formula $H^\alpha(M,\nu') = (M,\nu'\circ\alpha_M)$.
Conversely, if a functor
$H : {\mathcal M}^{\bf T'}\to{\mathcal M}^{\bf T}$
is such that $UH = U'$, where $U, U'$ are forgetful
and $F,F'$ are free $T$-algebra functors, then
$H$ induces a natural
transformation $\alpha^H : T \Rightarrow T'$
given by the composition
\begin{equation}\label{eq:alphaH}
T \stackrel{T\eta'}\longrightarrow TT' = UFU' F' = UFUHF'
\stackrel{U\epsilon HF'}\longrightarrow UHF' = U'F' = T'.
\end{equation}
These two rules are mutual inverses.

\nxpoint \label{s:mapMonads} More generally, 
given monad $\bf S$ in category $\MM$ and
monad $\bf T$ in category $\NN$,
a {\bf map of monads} $(K, \alpha):{\bf T}\to {\bf S}$
is a pair of a functor $K : \MM\to\NN$ and natural transformation
$\alpha : TK\tto KS : \MM\to\NN$ such that
$$\xymatrix{
TTK \ar[d]^{\mu^T K}\ar[r]^{T\alpha} & TKS \ar[r]^{\alpha S}
& KSS \ar[d]^{K\mu^S}\\
TK\ar[rr]^\alpha && KS
}$$
commutes and $\alpha\circ \eta^T K = K\eta^S : K\to KS$.
In $\CC\act$, the monads are now pairs 
$\tilde{\bf S} = ({\bf S},l^S)$,$\tilde{\bf T} = ({\bf T},l^T)$ 
and $K$ is replaced by a 
colax $\CC$-equivariant functor $(K,\zeta^K):\MM\tto\NN$, 
i.e. $\zeta^K_{C,M} : K(C\lozenge^\MM M)\tto C\lozenge^\NN KM$ 
form a binatural transformation of functors satisfying 
the coherences of types~(\ref{eq:colaxu}),(\ref{eq:colaxPsi}).

\nxpoint \label{s:mapPairs} 
A map of monads $(K,\alpha)$ is a {\bf map of pairs}
$(K,\alpha) : ({\bf T}, l^T)\to ({\bf S}, l^S)$ if
the following hexagon commutes
\begin{equation}\label{eq:pairsDeriv}\xymatrix{
TKG^\MM \ar[d]_{\alpha G^\MM}
\ar[r]^{T\zeta^K} & TG^\NN K \ar[r]^{l^T K} & G^\NN TK \ar[d]^{G^\NN \alpha} \\
KS G^\MM \ar[r]^{Kl^S} & KG^\MM S \ar[r]^{\zeta^K S} & G^\NN KS
}\end{equation}
where $G^\MM = (C \lozenge^\MM\_) \in\End_\CC(\MM)$ etc. (for all $C$). 

\nxpoint \label{s:compMapsMon} If $(K,\alpha) : {\bf T}\to {\bf S}$ and
$(L,\beta) : {\bf V}\to {\bf T}$ are two maps of monads,
then their composition is $(K,\alpha)\circ(L,\beta)  :=
(L\circ K, L\alpha \circ \beta K ) : {\bf V}\to {\bf S}$
which is again a map of monads as it follows by simple pasting:
$$\xymatrix{
VV(LK)\ar[rd]^{V\beta K} \ar[dd]^{\mu^V LK}
\ar[rr]^{V(L\alpha\circ\beta K)} &&
VLKS\ar[rr]^{(L\alpha\circ\beta K)S}\ar[rd]^{\beta KS}
&& LKSS \ar[dd]^{LK\mu^S}\\ &
VLTK \ar[r]^{\beta TK}\ar[ru]^{VL\alpha}
& LTTK\ar[d]^{L\mu^T K}\ar[r]^{LT\alpha} &
LTKS\ar[ru]^{L\alpha S} & \\
V(LK)\ar[rr]^{\beta K} && LTK\ar[rr]^{L\alpha} &&(LK)S
}$$

\nxpoint \label{s:compPairsMon}
For the equivariant case, there is nothing more here to show, as this 
makes sense in any 2-category. The composition of maps of pairs 
is in detail
$$
(L,\zeta^L,\alpha)\circ(K,\zeta^K,\beta) = 
(L\circ K, \zeta^L K\circ L\zeta^K, L\alpha\circ\beta K).
$$
The diagram expressing the fact that 
$$
L\alpha\circ\beta K : (VLK,l^VLK\circ V\zeta^L K \circ VL\zeta^K)
\tto (LKS, \zeta^L KS\circ L\zeta^K S\circ LKl^S)
$$
is $\CC$-equivariant may be obtained as follows:
$$\xymatrix{
VLKG^\MM\ar[d]_{\beta K G^\MM}\ar[r]^{VL\zeta^K} &
VLG^\NN K\ar[d]^{\beta G^\NN K}\ar[r]^{V\zeta^L K} &
VG^\PP LK\ar[r]^{l^V LK}& G^\PP VLK\ar[d]^{G^\PP \beta K}\\
LTKG^\MM\ar[d]_{L\alpha G^\MM}\ar[r]^{LT\zeta^K} &
LTG^\NN K\ar[r]^{Ll^T K} & 
LG^\NN TK\ar[d]_{LG^\NN\alpha}\ar[r]^{\zeta^L TK} &
G^\PP LTK\ar[d]^{G^\PP L\alpha}\\
LKSG^\MM\ar[r]_{LKl^S}& LKG^\MM S\ar[r]_{L\zeta^K S}&
LG^\NN KS\ar[r]_{\zeta^L KS}& G^\PP LKS
}$$

\nxpoint \label{s:transfMon} A transformation of maps of (usual) monads
$\sigma : (K,\alpha)\tto (L,\beta) : {\bf T}\to {\bf S}$
is a natural transformation $\sigma : K\tto L$
such that
\begin{equation}\label{eq:transfmapmon1}\xymatrix{
TK \ar[d]_{\alpha}\ar[r]^{T\sigma} & TL\ar[d]^{\beta} \\
KS\ar[r]^{\sigma S} & LS
}\end{equation}
commutes. For monads in $\CC\act$, te requirements ar the same, 
but of course the components need to be $\CC$-equivariant.
Thus the transformation of maps of pairs
$$
\sigma : (K,\zeta^K,\alpha)\tto(L,\zeta^L,\beta) 
: ({\bf T},l^T)\to({\bf S},l^S)
$$
is a usual transformation $\sigma : K\tto L$ 
satisfying the same Eq.~(\ref{eq:transfmapmon1}), 
but viewed as a transformation of pairs
$\sigma : (K,\zeta^K)\tto(L,\zeta^L) : ({\bf T},l^T)\to({\bf S},l^S)$,
is required to be $\CC$-equivariant, i.e. the square
\begin{equation}\label{eq:transfPairsEquiv}
\xymatrix{ 
K(C\lozenge^\MM M) \ar[r]^{\zeta^K_{C,M}}\ar[d]_{\sigma_{C\lozenge M}}
& C\lozenge^\NN KM\ar[d]^{C\lozenge \sigma_M}\\
L(C\lozenge^\MM M)\ar[r]_{\zeta^L_{C,M}} & C\lozenge^\NN LM
}\end{equation}
commutes for all $C$ in $\CC$ and $M$ in $\MM$.

\nxpoint (The cube for transformations of maps of pairs)
Denoting again, $G^\MM = C\lozenge^MM\_$ etc. we have the commutative
diagram
$$\xymatrix{
TKG^\MM\ar[r]^{T\zeta^K}\ar[d]_{T\sigma G^\MM} &
TG^\NN K\ar[d]^{TG^\NN\sigma}\ar[r]^{l^T K} &
G^\NN TK\ar[d]^{G^\NN T\sigma}\\
TLG^\MM\ar[r]_{T\zeta^L}&
TG^\NN L\ar[r]_{l^T L}&G^\NN TL
}$$
which is actually the upper face of the cube
\begin{equation}\label{eq:transCCact}\xymatrix{
TKG^\MM \ar[rr]^{T\zeta^K}
\ar[dr]^{T\sigma G^\MM}\ar[dd]^{\alpha G^\MM}
&& TG^\MM K\ar[rr]^{l^T K}\ar[rd]^{TG^\MM\sigma}
&& G^\NN TK\ar[dr]^{G^\NN T\sigma}\ar[dd]^{G^\NN\alpha} &
\\
& TLG^\MM\ar[rr]^{T\zeta^L}\ar[dd]^{\beta G^\MM} &&
TG^\MM L\ar[rr]^{l^T L} && G^\NN TL\ar[dd]^{G^\NN\beta}
\\
KSG^\MM\ar[rr]^{Kl^S} \ar[dr]_{\sigma S G^\MM} &&
KG^\MM S\ar[rr]^{\zeta^KS}\ar[rd]^{\sigma G^\MM S} 
&& G^\NN K S \ar[dr]^{G^\NN \sigma S} &
\\
& L SG^\MM \ar[rr]^{L l^S} && 
LG^\MM S\ar[rr]^{\zeta^L S} && G^\NN L S
}\end{equation}
where the bottom face is analogous commutative diagram
involving $L$ instead of $K$, where the left and right faces
commute because $\sigma$ is a transformation of usual monads,
and the front and back hexagons commute because $\beta$ and $\alpha$ are
maps of pairs, cf. diagram~(\ref{eq:pairsDeriv}). Hence the cube commutes.


\nxpoint \label{s:mixedpentagon} {\bf Theorem. (Mixed heptagon formula,
given a map of distributive laws) }
{\it
Let $l^S, l^T$ be two distributive laws between a $\CC$-actions and 
monads, ${\bf S}, {\bf T}$ in $\CC$-actegories $\MM, \NN$ respectively, 
and $(K,\zeta^K,\alpha) : ({\bf T}, l^T)\tto ({\bf S}, l^S)$ 
a map of pairs as above.
Then for each $C$ in $\CC$, $G^\MM := C\lozenge^\MM\_ \in \End_\CC(\MM)$,
the following diagram commutes
}
$$\xymatrix{
TKSG^\MM\ar[r]^{TKl^S} \ar[d]_{(K\mu^S \circ \alpha S) G^\MM}
& TKG^\MM S \ar[r]^{T\zeta^K S}
& TG^\NN KS
\ar[r]^{l^T K} & G^\NN TKS \ar[d]^{G^\NN(K\mu^S \circ\alpha S)}
\\
KSG^\MM\ar[rr]^{Kl^S} &&  KG^\MM S \ar[r]^{T\zeta^K} & G^\NN KS
}$$
{\it Proof.} This is obtained by the pasting of the following diagram
$$\xymatrix{
TKSG^\MM\ar[r]^{TKl^S} \ar[d]_{\alpha S G^\MM}
& TKG^\MM S \ar[r]^{T\zeta^K S} \ar[d]_{\alpha G^\MM S}
& TG^\NN KS
\ar[r]^{l^T K} & G^\NN TKS \ar[d]^{G^\NN\alpha S}
\\
KSSG^\MM \ar[r]^{KSl^S} \ar[d]_{K\mu^S G^\MM}
& KSG^\MM S \ar[r]^{K l^S S} &
KG^\MM SS \ar[r]^{\zeta^K SS} \ar[d]^{KG^\MM\mu^S}
& G^\NN K SS \ar[d]^{G^\NN K\mu^S}
\\
KSG^\MM\ar[rr]^{Kl^S} &&  KG^\MM S \ar[r]^{\zeta^K S} & G^\NN KS
}$$
where the upper left corner is commutative by naturality of $\alpha$,
the upper right by the pair property of $\alpha$,
the left lower corner is the pentagon for the distributive law $l^S$
and the right lower corner comes from the naturality of
$\mu^S$. Q.E.D.

\nxpoint Recall that a PRO is a strict monoidal category, whose object
part is the set of natural numbers (including 0) 
and the tensor product of objects is the addition of natural numbers
(and the unit object is $0$).
Different PRO-s differ by the morphisms,
and the tensor product on morphisms is still usually denoted by $+$
but typically it is not commutative. A (strict) representation of
PRO $\mathcal D$ in a monoidal category $\mathcal E$
is a strict monoidal functor $\DD\to\mathcal E$. There is
an obvious way to define PRO-s by morphism generators (under composition 
and ``addition'') and relations. 

We saw above that an endocell in $\CC\act$ is an endofunctor $T$ together
with a ``distributive law'' between $\CC$ and $T$ what is a
binatural transformation $l^T$ satisfying 
two commutative diagrams~(\ref{eq:colaxu}), (\ref{eq:colaxPsi}) with 
$F = T$ and $l^T = \zeta^F$. Given a representation
$\bm T^\bullet : \DD\to\mathcal E$
we denote by $\bm T^n := T(n)$ and simply $\bm T := \bm T(1)$

{\bf Theorem.} {\it 
A (strict) representation of a PRO 
$\tilde{\bf T}^\bullet : \DD\to\End_\CC(\MM)$ is the same
as a pair $({\bf T},l)$ where
${\bf T}^\bullet : \DD\to\End(\MM)$ is a representation 
and $l = l^T$ is a binatural transformation
$$
l^T : T(\_\lozenge\_)\tto \_\lozenge T(\_),
\,\,\,\,\,\, l^T_{C,M} : T(C\lozenge M)
\tto  C\lozenge T(M) 
$$
satisfying~(\ref{eq:colaxu}), (\ref{eq:colaxPsi})
and such that for every $\alpha : n\to m$ the (n+m+2)-gon
$$\xymatrix{
T^n (C\lozenge M)\ar[r]^{T^{n-1} l_{C,M}}
\ar[d]^{\alpha_{C\lozenge M}}
& T^{n-1} (C\lozenge TM)\ar[r]^{T^{n-2} l_{C,TM}}
& \ldots T(C\lozenge T^{n-1}M)\ar[r]^{l_{C,T^{n-1}M}}
& C\lozenge T^n M\ar[d]^{C\lozenge \alpha_M} \ar[d]^{C\lozenge \alpha_M}
\\
T^{m} (C\lozenge M)\ar[r]_{T^{m-1} l_{C,M}}
& T^{m-1} (C\lozenge TM)\ar[r]_{T^{m-2} l_{C,TM}}
& \ldots T(C\lozenge T^{m-1}M)\ar[r]_{l_{C,T^{m-1}M}}
& C\lozenge T^m M 
}$$
commutes.
}

The last condition simply says that 
$\alpha : T^n\tto T^m$ is in fact a
$\CC$-equivariant transformation $\alpha : (T^n,l^{T^n})\tto(T^{m}, l^{T^m})$ 
of colax $\CC$-equivariant endofunctors,
where $l^{T^n} := l T^{n-1}\circ \ldots \circ
T^{n-2} l T \circ T^{n-1} l$. 
This gives as many new diagrams
as there are many primitive natural transformations in the game.
For example a nonunital comonad has a coproduct $\delta$ hence
the distributive laws between $\CC$-action and a nonunital comonad
satisfy one more axiom, what amounts to 3 diagrams total. 
More precisely, one has a structure of a PRO on natural numbers where
$\delta$ etc. are the maps between $n$ and $m$ instead of $T^n$ and
$T^m$ and we are dealing in fact with a strict monoidal functor
from this PRO to the category of endofunctors of $\MM$ (called
also a strict representation of this PRO). Now I claim that a 
strict representation in $\End_\CC(\MM)$ is simply a pair
of a representation in $\End(\MM)$ and a distributive
law in the generalized sense, satisfying $n+2$ relations
if the PRO is generated by $n$ morphisms. 

\nxpoint Now specialize $\CC$ to the image of a representation 
$G_\bullet : \CC_0 \to \End(\MM)$ of (another) 
PRO $\CC_0$ in $\End(\MM)$. The generating object is $G_1 = G_\bullet(1)$.
$\CC$ is itself not necessarily a PRO even in this case, as
there may be a nonzero kernel of $G_\bullet$ on the level of objects,
but this presents no difficulty in the following.
This is a strict monoidal subcategory 
of $\End(\MM)$. 
Thus we have now two PRO-s in the game.
First of all in this case the diagrams~(\ref{eq:colaxu}), (\ref{eq:colaxPsi})
may be skipped all together! 
Namely $\Psi, u, l_{\bm 1,M}$ are all identities,
hence ~(\ref{eq:colaxu}) is a tautology, while~(\ref{eq:colaxPsi})
for general $A = G^n$, $B = G^m$ says simply
\begin{equation}\label{eq:lTargnm}
l_{G^{n+m}, M} = G^n (l_{G^m, M})\circ l_{G^n, G^m M}. 
\end{equation}
and in particular
\begin{equation}\label{eq:lTargn1}
l_{G^{n}, M} = G^{n-1} (l_{G, M})\circ l_{G^{n-1}, GM}. 
\end{equation}
what can be iterated to obtain  
\begin{equation}
l_{G^n, M} = G^{n-1} (l_{G,M})\circ G^{n-2}(l_{G, GM})\circ 
\ldots\circ G(l_{G,G^{n-2}M})\circ l_{G,G^{n-1} M}.  
\end{equation}
Thus {\bf every $l_{G^n,M}$ can be in the case when $\Psi$-s are
strict described in terms of $l_{G,G^s M}$ for varying $s\leq n$}. 
In particular, it is enough to consider the distributive laws
with one naturality  
$$
l : TG\tto GT, \,\,\,\,\,l_M : = l_{G,M}.
$$
We denote by $l^{(n)}_M := l_{G^n, M}$.
This way we have 
\begin{equation}\label{eq:ln}
l^{(n)} = G^{n-1} l \circ G^{n-2}lG\circ\ldots\circ GlG^{n-2}\circ lG^{n-1}.
\end{equation}
The naturality of $l_{C,M}$ in first argument, for
$\delta : G^n\to G^m \in\mathrm{Mor}\,\CC = G_\bullet(\CC_0)$ 
says that $(n+m+2)$-gon
\begin{equation}\label{eq:lGnat}\xymatrix{
TG^n\ar[d]_{T\delta} \ar[r]^{lG^{n-1}} &
GTG^{n-1} \ar[r]^{GlG^{n-2}} &
\ldots G^{n-1}TG \ar[r]^{\,\,\,\,\,\,\,G^{n-1}l} & 
G^n T\ar[d]^{\delta T}\\
TG^m\ar[r]^{lG^{m-1}} &
GTG^{m-1} \ar[r]^{GlG^{m-2}} &
\ldots G^{m-1}TG \ar[r]^{\,\,\,\,\,\,\,G^{m-1}l} & 
G^m T
}\end{equation}
commutes. 

From now on, whenever we discuss the distributive law between
two representations of PRO-s we will consider just the 
transformation $l$ with one naturality. 

For example, let $\CC$ be the PRO for counital coalgebras. Its 
set of morphisms is generated by
a morphism $\delta : 1\to 2$, satisfying the 
coassociativity $(\delta +\id)\delta = (\id+\delta)\delta$ 
and a morphism $\epsilon : 1\to 0$ satisfying $(\epsilon + \id)\delta = 
(\id+\epsilon)\delta = \id$. An action of this PRO is, of course, a 
counital comonad. Then,~(\ref{eq:lGnat}) becomes a pentagon for
$\delta$ and a triangle for $\epsilon$. 

More generally, if we have two endofunctors 
first with a structure arising from a representation of
one PRO and another with a structure arising from another PRO, 
with $k$ and $p$ relations respectively,
then we get in total $k+p$ additional diagrams for $l$ (there are
no conditions on $l$ except to be a transformation $TG\to GT$ otherwise).
The sizes of diagrams are always $n+m+2$ where $n$ and $m$ are 
the domain and codomain of a morphism in one or another PRO in question.

\nxpoint If $G = T$ is an underlying functor or a comonad $\bm G$ and
the distributive law $l : GG \tto GG$ satisfies 
the quantum Yang-Baxter equation $Gl\circ lG\circ Gl = lG\circ Gl\circ lG$
we say that $l$ is a strong braiding on the comonad $\bm G$. 
Then formula~(\ref{eq:ln}) defines a distributive law between 
$G$ and $G^{n}$ where the latter is inductively
equipped with a composite comonad structure using $l^{(p)}$.
for $p < n$. These results are discussed 
in our earlier article~\cite{skoda:cyclic}.

\nxpoint \label{s:mapsPairs}
Suppose $\tilde{\bm S}_\bullet :\PP\to\End_\CC(\MM)$, 
$\tilde{\bm T}_\bullet :\PP\to\End_\CC(\NN)$ are representations of
a fixed PRO $\PP$. As before, $\tilde{\bm S} = (\bm S_\bullet, l^S)$
and $\tilde{\bm T} = (\bm T_\bullet, l^T)$. 
A (colax) {\bf map of pairs} 
$(K,\zeta^K, \alpha) : ({\bf T}_\bullet, l^T)\to ({\bf S}_\bullet, l^S)$ 
is a colax $\CC$-equivariant functor 
$(K,\zeta^K) : \MM\to\NN$ together with a binatural transformation
$\alpha : TK\tto KS$ such that hexagon~(\ref{eq:pairsDeriv}) commutes
and such that for every morphism $\tau : n\to p$ in $\PP$
with $\tau^T := \bm T_\bullet(\tau)$ the following diagram
also commutes:
\begin{equation}\label{eq:multigonMapPairs}\xymatrix{
T^n K \ar[d]^{\tau^T K}
\ar[r]^{T^{n-1}\alpha} & T^{n-1} KS\ar[r]^{T^{n-2}\alpha S} &
\ldots TKS^{n-1}\ar[r]^{\,\,\,\,\,\,\,\alpha S^{n-1}} & KS^n\ar[d]^{K\tau^S} \\
T^m  K 
\ar[r]^{T^{p-1}\alpha} & T^{p-1} KS\ar[r]^{T^{p-2}\alpha S} &
\ldots TKS^{p-1}\ar[r]^{\,\,\,\,\,\,\,\alpha S^{p-1}} & KS^p
}\end{equation}
A map of pairs may be thought of as a colax $\CC$-equivariant
intertwiner from $\tilde{\bm S}_\bullet$ to $\tilde{\bm T}_\bullet$.

\nxpoint Generalizing the notation from~(\ref{eq:ln}) for
any natural transformation $\alpha  : TK\tto KS$ define
$\alpha^{(n)} := T^{n-1}\alpha \circ T^{n-2}\alpha S\circ \ldots\circ
\alpha S^{n-1} : T^{n}K\to KS^{n}$. Let $L : \NN\to\RR$ be a functor,
${\bf V}_\bullet : \PP\to\uEnd_\CC(\RR)$ 
a $\CC$-equivariant representation of $\PP$, and
$(L,\zeta^L,\beta) :
({\bf V}, l^V)\to({\bf T},l^T)$ a map of pairs. 

{\bf Lemma.} 
$L\alpha^{(n)} \circ \beta^{(n)} K = (L\alpha \circ \beta K)^{(n)}$.

This follows by easy induction. Using this one easily proves that 
the analogue of the multigon~(\ref{eq:multigonMapPairs})
for $L\alpha\circ\beta K$ is commutative. This together with 
\refpoint{s:compMapsMon} gives 

{\bf Proposition.} {\it
The rule
$$
(L,\zeta^L,\beta)\circ(K,\zeta^K,\alpha) :=
(L\circ K, \zeta^L K \circ L\zeta^K, L\alpha\circ\beta K) 
: ({\bf V},l^V)\to({\bf S},l^S)
$$
gives a (associative) composition 
of maps of pairs.
}

\nxpoint {\bf Theorem. (Mixed heptagon for maps of endofunctor 
$\CC$-equivariant representations of a PRO)}
{\it For every $(\tau : n\to p)\in\mathrm{Mor}\,\PP$, 
the following diagram 
$$\xymatrix{
TKS^{n-1} G^\MM \ar[d]^{(K\tau^S\circ \alpha S^{n-1}) G^\MM}
\ar@<2pt>[r]^{TKl^{S(n-1)}} &
TKG^\MM S^{n-1} \ar@<2pt>[r]^{T\zeta S^{n-1}} &
TG^\NN KS^{n-1} \ar@<2pt>[r]^{l^T KS^{n-1}} &
G^\NN TKS^{n-1}\ar[d]_{G(K\tau^S\circ \alpha S^{n-1})}
\\
KS^{p}G^\MM \ar[r]^{Kl^{S(p-1)}} & KG^\MM S^p\ar[rr]^{\zeta S^p} && G^\NN KS^p
}$$
commutes, where $l^{S(n-1)} = S^{n-2} l^S \circ \ldots \circ l^S S^{n-2}$
and $\zeta = \zeta^K$.
}

This proof is completely analogous 
to the proof of~\refpoint{s:mixedpentagon}, 
hence it is left to the reader. We call these identities
``mixed'' because unlike the diagrams for $l^T$ and $l^S$ separately,
they involve both $l^T$ and $l^S$. 

\nxpoint \label{s:transfpairs} 
Finally, the notion of transformation of maps of pairs
$\sigma : (K,\zeta^K,\alpha)\tto(L,\zeta^l,\beta) 
:(\bm T_\bullet, l^T)\to(\bm S_\bullet, l^S)$ is identical as 
in the case of monads in \refpoint{s:transfMon} 
(as it does not involve morphisms in $\PP$): 
require the commutativity of 
~(\ref{eq:transfmapmon1}) and~(\ref{eq:transfPairsEquiv}). 

\nxpoint {\bf Theorem. } {\it
The $\CC$-equivariant endofunctor representations of PRO $\PP$ 
in varying $\CC$-actegories
are objects of a 2-category $\mathrm{Rep}_{\CC\act}(\PP)$
where 1-cells are (colax) maps of pairs in the sense of \refpoint{s:mapsPairs}
and 2-cells are transformations of maps of pairs in the sense of
\refpoint{s:transfpairs}. We also consider the 2-subcategory
$\mathrm{Rep}_{\CC\actp}(\PP)\subset\mathrm{Rep}_{\CC\act}(\PP)$ 
where the 1-cells are those maps $(K,\zeta^K)$
of pairs whose coherences $\zeta^K$ are invertible.
}

The details are left to the reader.

\nxpoint (The category $\distr(\MM, G)$ of distributive laws
between an endofunctor (resp. a (co)monad) $G$ and
varying monads in a fixed category $\mathcal M$.) 
Objects of $\distr(\MM, G)$ are pairs
$({\bf T}, l)$ where $\bf T$ is a monad in $\mathcal M$ and
$l$ is a distributive law from $G$ to $T$.
Morphisms $({\bf T},l)\to ({\bf T},l')$
are the monad morphisms $\alpha : {\bf T}\to {\bf T'}$
such that there is the following commuting square of natural
transformations of endofunctors:
\begin{equation}\label{eq:distrMor}\xymatrix{
TG \ar[d]^{\alpha G}
\ar[r]^{l} & GT \ar[d]^{G \alpha}
\\
T'G\ar[r]^{l'} & GT'
}\end{equation}
It is clear that if we $\CC$ is the PRO with only trivial morphisms, 
and $G = \bm G(\bm 1)$ for a representation $\bm G :\CC\to\uEnd(\PP)$
then $\distr(\MM,G)$ is simply the full sub-1-category of 
(the decategorification of) $\mathrm{Rep}_{\CC\act}(\PP)$ whose 0-cells are
equivariant representations $\bm T$ of the PRO $\PP$ for monoids 
(i.e. monads) in $\MM$.

\nxpoint \label{s:Beck0} {\bf The original theorem of Beck.}

(i) Let $l : TG\tto GT$ be a distributive law from an endofunctor
(resp. monad) $G$ to a monad ${\bf T} = (T,\mu,\eta)$. Then the rule
\begin{equation}\label{eq:distr2lift}
\tilde{G}: (M,\nu)\mapsto(GM,\nu_l) =  (GM,G(\nu)\circ l_M),
\,\,\,\nu_l : TGM \stackrel{l_M}\to GTM \stackrel{G(\nu)}\to GM,
\end{equation}
defines an endofunctor on $\MM^{\bf T}$ lifting $G$ to an
endofunctor (resp. monad).

(ii) Conversely, if $U : \MM^{\bf T}\to \MM$
is the forgetful functors
(forgetting the monad action: $(M,\nu)\mapsto M$), and
$\tilde{G} : \MM^{\bf T}\tto\MM^{\bf T}$ and endofunctor
such that $U \tilde{G} = G U$ then for any object $M$ in
$\MM$, the composition
\begin{equation}\label{eq:lift2distr}\xymatrix{
TGM \ar[rr]^{TG(\eta_M)}&&
TGTM \ar[rr]^{U(\epsilon_{\tilde{G}FM})}&& GTM
}\end{equation}
defines the $M$-component of a distributive law $TG\tto GT$.

(iii) These two rules are mutual inverses.
\nxpoint \label{sec:alpha2HGeq}
{\bf Proposition.} {\it Condition~(\ref{eq:distrMor})
ensures that the induced functor $H^\alpha$
among the Eilenberg-Moore categories
will be $G$-equivariant.
}

{\it Proof.} Let $\tilde{G}$ and $\tilde{G}'$ be the lifts of $G$ in
$\MM^{\bf T}$ and $\MM^{\bf T'}$ respectively.

For all $(M,\nu')$ in $\MM^{\bf T}$,
$$\begin{array}{lcl}
H^\alpha(\tilde{G}'(M,\nu')) &=&
H^\alpha(GM,G(\nu')\circ l_M)\\
&=& (GM, G(\nu')\circ l_M\circ \alpha_{GM})
\\
&\stackrel{(\ref{eq:distrMor})}= &
(GM, G(\nu')\circ
G(\alpha_M)  \circ l'_M)\\
&=& (GM, G(\nu'\circ \alpha_M) \circ l_M)\\
&=& \tilde{G}(M, \nu'\circ \alpha_M) \\
&=& \tilde{G}H^\alpha(M,\nu').
\end{array}$$

\nxpoint \label{sec:lemmaepsPepsietaP}
{\bf Lemma.} {\it Given any functor $H: \MM^{\bf T'}\to \MM^{\bf T}$
satisfying $UH = U'$ the following identity holds
}
\begin{equation}\label{eq:epsPepsietaP}
H\epsilon' \circ \epsilon HF'UH \circ F\eta'UH = \epsilon H
\end{equation}
{\it Proof.} This follows from the naturality square
$\epsilon H \circ FUH\epsilon' = H\epsilon'\circ\epsilon HF'U'$
and the adjunction triangle $UH\epsilon\circ \eta' UH = \id_{UH}$
for $F' \vdash U' = UH$:
$$\xymatrix{
&&FUH\ar[lld]_{F\eta'UH}
\ar[d]_{\id}
\ar[rrrdd]^{\epsilon H}
&&
\\
FUH F' UH\ar[rd]_{\epsilon HF'UH}
\ar[rr]_{FUH\epsilon'} &&
FUH\ar[rrrd]_{\epsilon H}&&
\\
& HF'UH\ar[rrrr]^{H\epsilon'} &&&& H
}$$

\nxpoint \label{sec:mixedPhard}
{\bf Theorem.} ({\bf Mixed pentagon formula, given a functor $H$})
{\it Let $l,l'$ be two distributive laws
from an endofunctor $G$ on $\MM$ to the monads
${\bf T}, {\bf T}'$ respectively.
If $H : \MM^{\bf T'}\to\MM^{\bf T}$ is
a functor such that $UH = U'$ and $\tilde{G}H = \tilde G'$
then
\begin{equation}
\label{eq:D1M}
\xymatrix{
TT'G \ar[r]^{Tl'}
\ar[d]_{U\epsilon HF'G} &
T G T' \ar[r]^{l T'} &
  GTT'
\ar[d]^{GU\epsilon HF'}\\
T'G\ar[rr]^{l'}&& GT'.
}
\tag{D1M}
\end{equation}
}
Notice that two different distributive laws
(twice $l'$ and once $l$) appear
in the formula and that this formula (D1M) reduces to (D1) if
${\bf T} = {\bf T}'$ and $l = l'$ (then $H = \id$ and $\mu = U\epsilon F$).

{\it Proof.} By~(\ref{eq:lift2distr}), axiom (D1M) will follow
by the commutativity of
$$\xymatrix{TT'G
\ar[r]^{TT'G\eta'}\ar[d]^{U\epsilon HF'G}&
TT'GT'\ar[rr]^{TU'\epsilon' \tilde{G}'F'}
\ar[d]^{U\epsilon HF'GT'}&&
TGT'\ar[r]^{l T'}
\ar[rd]^{U\epsilon GF'}
& GTT' \ar[d]^{G\mu}\\
T'G\ar[r]^{T'G\eta'}
& T'GT'\ar[rrr]^{U'\epsilon' GF'}
&&&GT'
}$$
Here the left-hand square and the middle rectangle
commute by the naturality of $U\epsilon$.
The corner triangle on the left will be expanded further
to prove its commutativity:
$$\xymatrix{
TGT'\ar[rr]^{TG\eta T'}
\ar[rrdd]_{TG\eta' T'}
\ar[ddd]_\id
&& TGTT'\ar[rr]^{U\epsilon \tilde{G}FT'}
\ar[dd]^{TG\alpha^H T'}\ar[rd]^{TGT\eta'T'}
&& GTT'\ar[ddd]^{GU\epsilon HF'}
\\
&&
& TGTT'T'\ar[ld]^{TGU\epsilon HF'T'} &
\\
&& TGT'T'\ar[dll]^{TG\mu'} & &
\\
TGT'\ar[rrrr]^{U\epsilon H\tilde{G}F'} &&&& GT'
}$$
The upper horizontal line is expanded using~(\ref{eq:lift2distr})
and noticing $TGTT'= UFGUFT'= UFU\tilde{G}FT'$; what is sent
by $U\epsilon\tilde{G}FT'$ into $U\tilde{G}FT' = GUFT' = GTT'$.
The commutativity of the 3 triangles on the left is evident:
for the leftmost
follows from $\mu' \circ \eta' T = \id$ and the functoriality of
$TG$; for the next triangle by the unit axiom for $\alpha^H$,
i.e. $\eta' = \alpha^H\circ\eta$;
and for the third triangle
by the definition~(\ref{eq:alphaH}) of $\alpha^H$.
To prove that the right-hand hexagon
also commutes it is sufficient to prove that the 3 sides on the
right compose to $UFUX$ where $X = \tilde{G}\epsilon HF'$.
Indeed, $UX = GU \epsilon HF'$
is the top-down morphism on the right and the hexagon readily
reduces to a naturality rectangle for $U\epsilon$.
The 3 arrows on the right in fact compose to
$$\begin{array}{lcl}
TG(\mu'\circ U\epsilon HF'T'\circ T\eta'T')
&=&
UFG(UH\epsilon'F'\circ U\epsilon HF'U'F'\circ UF\eta'U'F')
\\
&=&
UFU\tilde{G}(H\epsilon'\circ\epsilon HF'U'\circ F\eta'U')F'
\\
&=&
UFU\tilde{G}(H\epsilon'\circ\epsilon HF'U'\circ F\eta'U')F'
\end{array}$$
The RHS is evidently equal to $UFUX$ as required
if the expression in the brackets equals $\epsilon H$.
This is exactly the content of the previous
Lemma~\refpoint{sec:lemmaepsPepsietaP}, i.e. formula~(\ref{eq:epsPepsietaP}).

This finishes the proof of the ``mixed pentagon formula''.

\nxpoint \label{sec:HeqalphaD} {\bf Corollary.} {\it
If $H : \MM^{\bf T'}\to\MM^{\bf T}$ is a functor satisfying
$UH = U'$, equivariant in the sense
$H\tilde{G} = \tilde{G}'H$, then $\alpha = \alpha^H$ is
a morphism in $\distr(\MM,G)$, i.e. it satisfies~(\ref{eq:distrMor}).
}

{\it Proof.} The required commutativity of~(\ref{eq:distrMor}),
by the definition $\alpha^H := U\epsilon HF' \circ T\eta'$,
reduces to the commutativity of the
external part of the diagram
$$\xymatrix{
& T G \ar[rr]^{l}
\ar[ld]_{T\eta' G}
\ar[rd]^{TG\eta'}
&& GT \ar[rd]^{GT\eta'}
&
\\
TT'G\ar[rr]_{Tl'}
\ar[d]^{U\epsilon HF'G}&&
TGT'\ar[rr]_{lT'} &&
GTT'\ar[d]_{GU\epsilon HF'}
\\
GT'\ar[rrrr]^{l'} &&&& GT'
}$$
The commutativity of the left-top triangle
is the unit axiom for the distributive law $l$,
the right-top rectangle is commutative by the naturality of $l$,
and the bottom is the pentagon (D1M) from~\refpoint{sec:mixedPhard}. Q.E.D.

\nxpoint {\bf Theorem. (Mixed pentagon formula, given a map
$\alpha$ of distributive laws)} {\it Let $l, l'$ be two distributive
laws from an endofunctor $G$ to to monad ${\bf T, T'}$ respectively
and $\alpha : ({\bf T}, l) \tto ({\bf T}', l')$ a morphism
in $\distr(\MM,G)$. Then the following diagram commmutes
}
\begin{equation}
\label{eq:D1Ma}
\xymatrix{
TT'G \ar[r]^{Tl'}
\ar[d]_{(\mu'\circ\,\alpha T')G} &
T G T' \ar[r]^{l T'} &
  GTT'
\ar[d]^{G(\mu'\circ\,\alpha T')}\\
T'G\ar[rr]^{l'}&& GT'.
}
\tag{D1Ma}
\end{equation}
{\it Proof.} This simple proof is due M. Jibladze (personal communication).
Recall that $\mu = U\epsilon F$. Then the following
diagram is commutative:
$$\xymatrix{
TT'G \ar[r]^{Tl'}
\ar[d]_{\alpha T'G} &
T G T' \ar[r]^{l T'}\ar[d]^{\alpha GT'} &
  GTT'
\ar[d]^{G\alpha T'}\\
T'T' G\ar[d]\ar[r]^{T' l'}
\ar[d]_{\mu' G} &
T'GT' \ar[r]^{l' T'} & GT'T'\ar[d]^{G\mu'}
\\
T'G\ar[rr]^{l'}&& GT'
}$$
Indeed, the lower pentagon is a part of the statement that $l'$
is a distributive law. The left upper corner square is a
naturality square for $\alpha$. Finally the right upper corner is
expressing the condition that $\alpha$ is a map in $\distr(\MM,G)$
(composed by $T'$). The external part of this diagram evidently
gives~(\ref{eq:D1Ma}). Q.E.D.
\nxpoint {\bf Proposition.} {
If $H = H^\alpha$ in (D1M), or equivalently, by 2, $\alpha = \alpha^H$,
then the vertical arrows in (D1M)
are identical to the corresponding compositions of vertical arrows
in (D1Ma).
}

{\it Proof.}  It is sufficient to show
$U\epsilon H  F' = \mu' \circ \alpha^H T'$
as this implies the assertion both for the left-hand and right-hand
vertical arrows. In fact we show the stronger assertion that
$U\epsilon H = \mu' \circ \alpha^H U'$. Setting
$\mu' = U'\epsilon' F' = UH\epsilon' F'$
and $\alpha^H = U\epsilon HF' \circ UF\eta'$
we reduce the required identity to
$U\epsilon H = UH\epsilon' F'\circ U\epsilon HF' U'\circ UF\eta'U'$.
By the functoriality of $U$, the assertion follows from
Lemma~\refpoint{sec:lemmaepsPepsietaP},
that is formula~(\ref{eq:epsPepsietaP}). Q.E.D.

\nxpoint {\bf Theorem.} {\it
Given an endofunctor (resp.~comonad) $G$ in a category $\MM$,
the category $\distr(\MM, G)$ is canonically isomorphic to the category of
Eilenberg-Moore categories of varying monads equipped with
a lift of $G$, and functors commuting with the forgetful functors
and intertwining the lifts of $G$.
}

{\it Proof.} {\bf (i)} ({\bf Bijection for objects}) By the definition,
Eilenberg-Moore categories of $\bf T$-modules $\MM^{\bf T}$
are trivially in 1-1 correspondence with the monads ${\bf T}$,
and for a fixed monad the distributive laws are in bijection
with lifts by Beck's theorem \refpoint{s:Beck0}.

{\bf (ii)} ({\bf Bijection of $\mathrm{Hom}$-sets})
Given a pair of monads ${\bf T}$, ${\bf T}'$,
it is also classical that morphism of monads
are in 1-1 correspondence $\alpha \mapsto H^\alpha$
with the functors of Eilenberg-Moore categories
commuting with the forgetful functor.
So to show the bijection for morphisms there is only
one nontrivial thing to prove: the property that a map $\alpha$ of monads
is actually a morphism in $\distr(\MM,G)$ corresponds exactly to the
fact that $H^\alpha$ is intertwining the corresponding lifts of $G$.
But all the hard work there has been already done:
Proposition~\refpoint{sec:alpha2HGeq}, states this in one direction,
and Corollary~\refpoint{sec:HeqalphaD} does the converse.

{\bf (iii)} ({\bf $\alpha\mapsto H^\alpha$
is a contravariant functor})
This is certainly known, but we do not know the reference.
First of all, the identity functor $H = \id$ gives
$\alpha^\id = \id$ as it is clear by the adjunction triangle
$\epsilon F \circ F\eta$.
In the situation
$$
{\mathcal M}^{\bf T''} \stackrel{H'}\longrightarrow
{\mathcal M}^{\bf T'} \stackrel{H}\longrightarrow
{\mathcal M}^{\bf T}
$$
with $UH = U'$, $U'H' = U''$,
we need to show that $\alpha^{H'} \circ \alpha^{H} = \alpha^{H\circ H'}$.
The LHS is the composition
$$
UF\stackrel{UF\eta'}\longrightarrow
UFU'F'\stackrel{U\epsilon HF'}\longrightarrow
UHF'\stackrel{U'F'\eta''}\longrightarrow
U'F'U''F''\stackrel{U'\epsilon' H' F''}\longrightarrow
U'H'F'' = T''
$$
By naturality of $U\epsilon$ we may interchange
$U'F'\eta''\circ U\epsilon HF' = U\epsilon HF'U''F''\circ UFUHF'\eta''$
and furthermore interchange
$U'\epsilon' H' F'' \circ U\epsilon HF'U''F'
= U\epsilon HH'F''\circ UFU'\epsilon'H'F''
: UFU'F'U''F''\to U'H'F'' = T''$.
Thus we obtain that LHS equals
$$
UF\stackrel{UF\eta'}\longrightarrow UFU'F'
\stackrel{UFUHF'\eta''}\longrightarrow UFU'F'U''F''
\stackrel{UFU'\epsilon'H'F''}\longrightarrow UFU''F''
\stackrel{U\epsilon HH'F''}\longrightarrow T''
$$
Now the composition of the second and third morphism is $UF\alpha^{H'}$
by the definition, and $\alpha^{H'}\circ \eta' = \eta''$
hence the composition
of the first three transformations is
$UF\eta''$, therefore all 4 compose to the
$U\epsilon HH'F''\circ UF\eta'' = \alpha^{H\circ H'}$
by the definition. Q.E.D.

\nxpoint It is again standard that maps of monads $\alpha:{\bf T}\to{\bf S}$
are in 1-1 correspondence with the functors
$H : \MM^{\bf S}\to\NN^{\bf T}$, such that $U^T H = KU^S$.
$$\xymatrix{
\MM^{\bf S}\ar[d]_{U^S} \ar[r]^H & \NN^{\bf T} \ar[d]^{U^T}\\
\MM \ar[r]^K &\NN
}$$
We will below need the explicit formulas for this bijection.
Given a functor $H$ as above, the
corresponding map of monads $\alpha^H : TK\tto KS$
(cf. Borceux, II 4.5.1)
is the composition
$$
TK \stackrel{TK\eta^S}\longrightarrow TKS = U^T F^T K U^S F^S
= U^T F^T U^T H F^S \stackrel{U^T\epsilon^T H F^S}\longrightarrow
U^T H F^S = KU^S F^S = KS
$$
Conversely, given a morphism of monads
$\alpha$ we obtain the lift $H^\alpha$ simply as
$$
H^\alpha(M,\nu) := (KM, K(\nu)\circ \alpha_M).
$$
All together this is a canonical bijection,
clearly extending the formulas in~\refpoint{sec:alphaHcorr}.
Moreover, for any transformation of maps of monads
$\sigma : (K,\alpha)\tto (K',\alpha') : {\bf T}\tto{\bf S}$
one defines a natural transformation
$\tilde\sigma : H^\alpha\tto H^{\alpha'}$ by
$\tilde\sigma = \sigma U^T$, i.e.
$$
\tilde\sigma_{(M,\nu)} := \sigma_M : (KM,K(\nu)\circ\alpha_M)
\to (K'M,K'(\nu)\circ\alpha_M).
$$
We leave for the reader to check that $\tilde\sigma_{(M,\nu)}$ is really
a morphism in $\NN^{\bf T}$, i.e.
$$
\sigma_M \circ K(\nu)\circ\alpha_M =
K'(\nu)\circ\alpha'_M\circ T(\sigma_M).
$$
This transformation lifts $\sigma$, when
considered just as a transformation of
functors $\sigma : K\tto K'$. That means $\tilde\sigma = \sigma U^T$.

Conversely, given any natural transformation
$\theta : H^\alpha\tto H^{\alpha'}$ such that
$U^T(\theta_{(M,\nu)}) : KM\to K'M$ does not depend on $\nu$
and hence lifts a (unique) transformation of functors
$\theta_* : K\tto K'$, then $\theta_*$ is automatically
given by formula $(\theta_*)_M = U^T(\theta_{(TM,\mu_M)})$
which is a transformation of maps of monads
$$
\theta_* : (K,\alpha)\tto(K',\alpha') : {\bf T}\to {\bf S}
$$
Now we claim that $H$ is equivariant (intertwines $G^\MM$ and
$G^\NN$)
iff $\alpha$ is a map of pairs, i.e.~(\ref{eq:pairsDeriv}) holds.
Of course, if the coherence $\zeta$ is non-trivial then 
one needs to equip also $H$ with a coherence.
Moreover, one can consider a certain 2-category of small
categories each equipped with an
endofunctor $G$, a monad, say ${\bf T}$, and a distributive law;
with the maps of pairs as morphisms and certain class of
compatible modifications of such morphisms. Then there is a
2-isomorphism with a 2-category of Eilenberg-Moore categories,
equipped with lifts, equivariant functors of such and  their equivariant
natural transformations where everything commutes with the forgetful
functors.

\nxpoint {\bf Proposition.} {\it Condition~(\ref{eq:pairsDeriv}) ensures
a 2-cell $H^\alpha G^\MM \Rightarrow G^\NN H^\alpha$.
}

{\it Proof.} For all $(M,\nu) \in \MM^{\bf S}$,
$$\begin{array}{lcl}
H^\alpha G^\MM (M,\nu) & = & H^\alpha (G^\MM M, G^\MM(\nu)\circ l^S_M)
\\
&=& (KG^\MM M, KG^\MM(\nu)\circ K(l^S_M)\circ \alpha_{GM})
\\
&=& (KG^\MM M, KG^\MM(\nu)\circ G^\NN(\alpha_M)\circ l^T_{KM})
\\
&\Rightarrow& (G^\NN K M, G^\NN (K(\nu)\circ\alpha_M)\circ l^T_{KM})
\\
&=& G^\NN(KM, K(\nu)\circ\alpha_M)
\\
&=& G^\NN H^\alpha (M,\nu)
\end{array}$$
We used in the middle step 
the 2-cell $\zeta^K_M : KG^\MM M\tto G^\NN KM$ in the first component
and composing with it in the second component. 


\nxpoint {\bf Theorem.} {\it
The natural transformation $\tilde\sigma : H^\alpha\tto H^{\alpha'}$
induced from a transformation of monads $\sigma$ is equivariant iff
$\sigma : \alpha\tto\alpha'$ is a transformation of maps of pairs.
}

\nxpoint {\bf Theorem.} {\it 
If $\PP$ is the PRO for monoids then
2-category $\mathrm{Rep}_{\CC\act}(\PP)$ is isomorphic to the
following 2-category: 
the objects are triples $(\MM,\bm T,U^T : \MM^T\to\MM)$
where $\bm T$ is a monad in a $\CC$-actegory $\MM$,
$\MM^{\bm T}$ is the Eilenberg-Moore category of $\bm T$ 
equipped with a $\CC$-action making $U^T$ a strict monoidal
functor; 1-cells are colax $\CC$-equivariant functors of
Eilenberg-Moore categories
$\MM^{\bm T}\to\NN^{\bm S}$ commuting with the forgetful functor
and 2-cells the natural transformations of 
colax $\CC$-equivariant functors.
}

\nxpoint In his classical article~\cite{street:formalmon} {\sc R. Street}
has considered monads and Eilenberg-Moore objects in general
2-categories. The fact that the Beck's bijection between
lifts and distributive laws extends to an isomorphism of 2-categories,
may be viewed, after applying our correspondence between the 2-category
of distributive laws and the 2-category of equivariant monads, as
the correspondence between the Eilenberg-Moore objects and monads
inside the 2-category $\CC\act$. For this one needs to apply a result
on the existence of Eilenberg-Moore objects in this setup. S. Lack has
proved a general result of this type, namely existence of certain lax
limits whose combinations include the Eilenberg-Moore objects, in the 
2-category of pseudoalgebras over a 2-monad. In our case the 2-monad is
a strictification of the pseudomonad $\CC\times$ on $\cat$, whose
structure is induced from the monoidal category structure on $\CC$,
and whose pseudocoalgebras are coherent $\CC$-actions. In a way this
is more general than our approach as it allows other 2-monads: on the 
other hand our case is more general as the monads are generalized to
actions of PRO-s and more general $\DD$-actions. Some subtleties of the
latter case are discussed in \cite{skoda:biact}. Our approach also
emphasizes on explicit formulas for all the correspondences and
isomorphisms instead of equivalences at certain places. 

\nxpoint (Relative distributive laws) 
Recall that a pseudomonad in a 
$\mathrm{Gray}$-category $\mathcal K$ is 
an object $H$ in $\mathcal K$ and 
a pseudomonoid in the $\mathrm{Gray}$-monoid $\mathcal{K}(H,H)$
(for $\mathrm{Gray}$-pseudomonoids see e.g.~\cite{DayStreetMonBicat}).
Thus a pseudomonad is a tuple 
$\bm D = (D,\mu,\eta,\alpha^l,\alpha^r,\alpha^\mu)$ where
$D : X\to X$ is a 1-cell in $\mathcal K$, $\mu :DD\to D$ and
$\eta : D\to DD$ are 2-cells in $\mathcal K$ and 
the coherence for right unit $\alpha^r : \mu\circ D\eta\tto\id_D$,
the coherence for left unit $\alpha^l : \mu\circ \eta D\tto\id_D$
and the coherence for associativity
$\alpha^\mu : \mu\circ(D\mu)\tto\mu\circ(\mu D)$
are invertible 2-cells in $\mathcal K$ satisfying 
2 standard coherence identities.
Suppose we are given pseudomonads $\bm C$ and $\bm D$ in $\mathcal K$,
and a {\it fixed} 1-cell $X$ in $\mathcal K(H',H)$, for some
object $H'$ in $\mathcal K$. Suppose that $X$ is both the
$\bm C$-pseudoalgebra $(X,\rho,\psi^C,\xi^C)$
and $\bm D$-pseudoalgebra $(X,\nu,\psi^D,\xi^D)$, one may ask what makes 
the $\bm D$-pseudoalgebra structure (say colax-) $\bm C$-equivariant 
in the sense
that the defining 1-cell $\nu : DX\to X$ and the invertible 2-cells 
$$\psi^D : \nu\circ\eta_X \tto \Id_X,
\,\,\,\,\,\,\,\,\,\,\chi^D : \nu\circ(D\nu)\tto \nu\circ\mu_X, $$
are equipped with a structure of 1-cell and 2-cells in 
the 2-category of $\bm C$-pseudoalgebras, colax morphisms of pseudoalgebras,
and their natural transformations. For this to make sense we need
also a $\bm C$-structure on $DX$ what may need another distributive
law, but in many cases this part of the data is in fact 
canonically provided, while the additional structure above is not.
For example, if the pseudomonads are the cartesian products with 
monoidal categories then we can just use the commutativity of the 
cartesian product to identify $DCX$ and $CDX$ while their
actions on concrete $X$ does {\it not} trivially commute and what we 
discuss here is precisely the additional distributive structure for
the two actions. More generally, we can consider just some
``higher'' distributive law between the pseudomonads, 
$\mathrm{can} : DC\to CD$ and define the distributive laws for
pseudoalgebras relatively to it. 
For the 1-cell $\nu : DX\to X$ the additional structure is a 2-cell 
$$\tau : \rho\circ C\nu\tto\nu\circ D(\rho)\circ \mathrm{can}$$ 
in $\mathcal K$, where two coherences hold for $\tau$, namely
$$\xymatrix{
& \ddrtwocell<\omit>{<0> \tau}
CDX\ar[r]^{C\nu}\ar[d]^{\mathrm{can}} 
& CX\ar[dd]^\rho\\
DX\ar[rd]_=\ar[ru]^{\eta^C_{DX}}
\ar[r]^{D(\eta^C)}\drtwocell<\omit>{<-2> \,\,\,\,\,\,\,\,\,D(\psi^C)} 
& DCX\ar[d]^{D(\rho)}&\\
& DX\ar[r]^\nu& X
}=
\xymatrix{DX\ar[r]^{\eta^C_{DX}}\ar[d]_{\nu}
\drtwocell<\omit>{<0> \eta^C_\nu}&
CDX\ar[d]^{C(\nu)}\\
X\ar[r]^{\eta^C_X}\ar[rd]_=
\drtwocell<\omit>{<-2>\,\,\,\,\,\psi^C} & CX\ar[d]^{\rho}\\
& X
}$$
and the pasting 
$$\xymatrix{
CCDX\ar[dd]^{\mu^C_{DX}}\ar[rd]^{\mathrm{can}}\ar[rr]^{CC\nu} &
\drtwocell<\omit>{<0> C\tau}
& CCX\ar[rrdd]^{C\rho} && \\
&CDCX\ar[rd]^{CD\rho}\ar[d]^{\mathrm{can}}&&&\\
CDX\ar[rd]^{\mathrm{can}}& DCCX\ar[d]^{D\mu^C_X}\ar[rd]^{DC\rho}
& CDX\ar[d]^{\mathrm{can}}\ar[rr]^{C\nu}\ddrrtwocell<\omit>{<0> \tau}&
& CX\ar[dd]^\rho \\
&\drtwocell<\omit>{<-2> \,\,\,\,\,\,\,\,\,\,\,D(\chi^C_{X})}
DCX\ar[rd]_{D\rho}& DCX\ar[d]^{D\rho} && \\
&& DX \ar[rr] && X 
}$$ equals the pasting
$$\xymatrix{
CCDX\ar[dd]_{\mu^C_{DX}}\ar[rr]^{CC\nu}
\ddrrtwocell<\omit>{<0>\,\,\,\,\,\,\,\,\,\,\,\,(\mu^C_\nu)^{-1}}&& 
CCX\ar[dd]_{\mu^C_X}\ar[rd]^{C\rho} &\\
&&\drtwocell<\omit>{<0> \chi^C}& CX\ar[dd]^{\rho}\\
CDX\ar[d]^{\mathrm{can}}\ar[rr]^{C\nu} 
&\drtwocell<\omit>{<0> \tau}& CX\ar[rd]^{\rho}&\\
DCX\ar[r]_{D\rho}& DX\ar[rr]_{\nu}&& X
}$$
These coherences say precisely that
$(\psi^D,\tau): (X,\rho,\psi^C,\chi^C)\to (X',\rho',\psi'^C,\chi'^C)$ 
is a colax morphism of $\bm C$-pseudoalgebras. Notice that 
if the pseudonaturality of $\mu^C$ and $\eta^C$ is in fact 
naturality then we exactly get one triangle and one pentagon for
the nonidentity 2-cells. 
For the 2-cells $\psi^D : \nu\circ\eta^D_X\tto\id_X$ and 
$\chi^D : \nu\circ D(\nu) \tto \nu\circ\mu^{D}_X$ 
there is no additional structure
but rather a {\it requirement} that 
they are natural transformations of colax $\bm C$-equivariant morphisms of
$\bm C$-pseudoalgebras, what boils down to a bit expanded tin-can diagrams:
$$\xymatrix{
CX\ar[r]\ar[dd]^\nu \ar@/_2pc/[rr]_{C(\id_X)}
\drrtwocell<\omit>{<0>\,\,\,\,\,\,\,\,\,\,\,C(\psi)}
& CDX\ar[r] & CX \ar[dd]^\nu\\
& 
& \\
X\ar@/_2pc/[rr]_{\id_X}&& X
}\xymatrix{ \\ = \\ }
\xymatrix{CX\ar[r]^{C(\eta^D_X)}\ar[ddd]^\rho
\ar[rd]^{\eta^D_{CX}}& CDX\ar[d]^{\mathrm{can}}\ddrtwocell<\omit>{<0>\tau}
 \ar[r]^{C(\nu)} 
& CX\ar[ddd]^\rho\\
\drtwocell<\omit>{<0>\,\,\,\eta^D_\rho}&DCX\ar[d]^{D(\rho)}&\\
\drrtwocell<\omit>{<0>\,\,\,\,\,\psi}
&DX\ar[rd]^\nu&
 \\
X\ar[ru]^{\eta^D_X}\ar[rr]_{\id_X}&&X }
$$
If $\eta^D$ is again natural (in particular $\eta^D_\rho = \id$),
this identity boils down to a triangle
for natural transformations. 
The tin can identity for $\chi^D$ is as follows
\begin{equation}\label{eq:chipasting}\xymatrix{
CDDX\ar[rd]^{C\mu^D_X}\ar[d]^{\mathrm{can}}\ar[r]^{CD\nu}
\drrtwocell<\omit>{<0>\,\,\,\,\,\,\,\,\,\,\,C(\chi^D)}
& CDX\ar[r]^{C\nu}
& CX\ar[ddd]^\rho \\
DCDX\ar[d]^{\mathrm{can}}& CDX\ddrtwocell<\omit>{<0>\tau} \ar[ru]_{C\nu}
\ar[d]^{\mathrm{can}}&\\
DDCX\drtwocell<\omit>{<0>\,\,\,\,\,\,\,\,\,(\mu^D_\rho)^{-1}}
\ar[d]_{DD\rho}\ar[r]^{\mu^D_{CX}}&DCX\ar[d]^\rho&\\
DDX\ar[r]_{\mu^D_X}&DX\ar[r]_\nu & X
}
\xymatrix{ \\ = \\ }
\xymatrix{ CDDX\ar[r]^{CD\nu}\ar[d]^{\mathrm{can}}&
CDX\ar[r]^{C\nu}\ar[d]^{\mathrm{can}}
&CX\ar[ddd]^\rho
\\DCDX\ar[d]^{\mathrm{can}}\ar[r]^{DC\nu}
\drtwocell<\omit>{<0>\,\,\,\,\,\,\,D\tau}&
DCX\drtwocell<\omit>{<-1>\tau}\ar[d]^{D\rho}&\\
DDCX\drrtwocell<\omit>{<0>\,\,\,\,\chi^D}
\ar[d]_{DD\rho}& DX\ar[dr]^\nu&\\
DDX
\ar[r]_{\mu^D_X}\ar[ru]^{D\nu}&DX\ar[r]_\nu&X\\
}\end{equation}
Again in the 2-categorical situation, when $\mu^D_\rho$ is the identity
this boils down to a pentagon for natural transformation.
The distributive laws between two
actions of monoidal categories on a fixed category $X$ 
are a special case of this construction. Notice that each of the two 
pentagons and two triangles, is defined using a pasting diagram
which contains embedded exactly one pentagon or triangle for the higher
distributive law. This is an interesting ``recursive'' structure. 
We see that the distributive laws
between the pseudoalgebras are defined
relative to a higher distributive law $\mathrm{can}:CD\to DC$ 
between their pseudomonads which is in our case
``canonical'' and invertible,
but it may be not so. Moreover, the higher distributive law
may be in fact a pseudodistributive law as in~\cite{marmolejo:pdistr1},
and we again, {\it mutatis mutandis}, define the distributive laws between
the pseudoalgebras using essentially the same ``relative''
pasting diagrams as above, sometimes
with nontrivial 2-cells inserted in place of trivial ones. 
For example, the upper left pentagon in the left-hand 
diagram in~(\ref{eq:chipasting}) 
is then filled with a nontrivial 2-cell. 

\footnotesize{
\vskip .3in
{\bf Acknowledgements}. 
The main results of this paper have been obtained at Newton Institute,
Cambridge, during my stay in October 2006.
I thank the institute for providing the good and stimulating working
conditions; Profs. S. Majid, A. Connes and A. Schwarz for
the invitation and Prof. M. Jibladze for his generous
mathematical advice, as usual. The article has been written up partly
also at MPI Bonn, IH\'ES and IRB Zagreb. 
I also thank DAAD-Croatia project 
``Nonabelian cohomology and applications''
for the travel expenses Zagreb-Bonn.

}
\end{document}